\documentclass[12pt,a4paper,twoside,final,notitlepage,leqno]{article}
\usepackage[english]{babel}
\usepackage[T1]{fontenc}
\usepackage{epsfig, graphicx, amssymb}
\usepackage{amsmath,amsthm,epsfig,amsfonts}
\usepackage{float}
\usepackage{color}
\usepackage{graphicx}
\usepackage{diffcoeff,xfrac,float,cite,booktabs,multirow,collcell}
\usepackage{siunitx}
\def\br {\break}

\newcommand{\monitem}{ \smallskip \noindent $\bullet$ \quad  }
\newcommand{\moneq}{\vspace*{-7pt} \begin{equation} \displaystyle }
\newcommand{\moneqstar}{\vspace*{-6pt} \begin{equation*} \displaystyle }
\newcommand{\monendstar}{\vspace*{-6pt} \end{equation*}   }
\newcommand{\monend}{\vspace*{-7pt} \end{equation}   }
\newcommand{\moneqarraystar}{ \begin{eqnarray*} \displaystyle }
\newcommand{\monendarraystar}{ \end{eqnarray*}   }

\newcommand{\dd}{{\rm d}}

\newcommand{\tran}{{^{\rm t}}}
\def\br {\break}

\definecolor{vertfonce}{rgb}{0.0, 0.5, 0.0}

\def\section*#1{}

\usepackage{fancyhdr}
\renewcommand{\headrulewidth}{0pt}

\begin{document}

\fancypagestyle{plain}{ \fancyfoot{} \renewcommand{\footrulewidth}{0pt}}
\fancypagestyle{plain}{ \fancyhead{} \renewcommand{\headrulewidth}{0pt}}

~

  \vskip 2.1 cm

\centerline {\bf \LARGE {A variational symplectic scheme}}

 \bigskip 

\centerline {\bf \LARGE {based on Lobatto's quadrature}}

\bigskip  \bigskip \bigskip

\centerline { \large   Fran\c{c}ois Dubois$^{ab}$ and Juan Antonio Rojas-Quintero$^{c}$}

\smallskip  \bigskip

\centerline { \it  \small
  $^a$ Laboratoire de Math\'ematiques d'Orsay, Facult\'e des Sciences d'Orsay,}

\centerline { \it  \small   Universit\'e Paris-Saclay, France.}

\centerline { \it  \small
$^b$    Conservatoire National des Arts et M\'etiers, LMSSC laboratory,  Paris, France.}

\centerline {\small
orcid  0000-0003-4858-6234, francois.dubois@lecnam.net} 

\centerline { \it  \small
  $^c$  Secihti/Tecnológico Nacional de México, I.T. Ensenada,}

\centerline { \it  \small
    22780 Ensenada, BC, Mexico.}

\centerline {\small
orcid  0000-0002-2311-6933, jarojas@secihti.mx}



\bigskip  \bigskip

\centerline {{28 October  2025}
{\footnote {\rm  \small $\,$ This contribution,  presented 
at the 7th International Conference on  {\it Geometric Science of Information},
Saint-Malo, 29-31 October 2025, is published in
{\it Lecture Notes in Computer Science}, volume 16034, pages  332-342, 2025.}}}




\bigskip  \bigskip
\noindent {\bf \large Abstract}

\noindent
  We present a variational integrator based on the Lobatto quadrature for the time
  integration of dynamical systems issued from the least action principle.
  This numerical method uses a cubic interpolation of the states and the action is approximated
  at each time step by Lobatto's formula. Numerical analysis is performed on {both a} harmonic oscillator 
  {and a nonlinear pendulum}.
  The {geometric} scheme is conditionally stable, sixth-order accurate, and symplectic.
  It preserves an approximate energy quantity. Simulation results illustrate the performance 
  {and the superconvergence} 
  of the proposed method. 

\noindent

\newpage

\bigskip \bigskip    \noindent {\bf \large    1) \quad  Introduction}
%
%

\smallskip 
A well-known geometric property of Hamiltonian systems is that their flows preserve the phase space volume,
a Poincaré invariant. However, integrable Hamiltonian systems are rare. Numerical methods are usually required
to find solutions. By generating functions of canonical transformations, symplectic variants
of standard integrators can be achieved.
We refer {\it e.g.} to the  symplectic Euler or Runge-Kutta \cite{Hairer2006}, 
or the implicit midpoint scheme that coincides with the average
constant acceleration Newmark variant  for linear systems. This method was proposed by J.C. Simo \cite{simo1992}
and is symplectic for both the linear and nonlinear cases.
However, a good way to embed the natural geometry of Hamiltonian systems into an integrator is to start
with the variational principle of least action.

\smallskip \noindent 
In this contribution, we first recall the Lobatto quadrature scheme in Section 2. 
Section~3 
presents the selected interpolation of functions with finite elements. Then, starting with the least action principle,
the Lobatto integrator is developed on the harmonic oscillator in Section 4. 
The resulting scheme
is expressed as a variant of the implicit midpoint integrator~\cite{simo1992} and the Simpson  
integrator proposed
in \cite{dubois2023_gsi}; it is a special case of the Galerkin methods introduced in \cite{marsden2001}.
The symplectic structure of the Lobatto integrator is analyzed at the end of Section 4. 
The method preserves a discrete Hamiltonian as remarked by \cite{benettin1994}.
Some numerical results are presented in Section 5 
before some conclusive words.

\bigskip \bigskip    \noindent {\bf \large    2) \quad  Lobatto's Quadrature}

\fancyhead[EC]{\sc{Fran\c{c}ois Dubois and Juan Antonio Rojas-Quintero}}
\fancyhead[OC]{\sc{A variational symplectic scheme based on Lobatto's quadrature}}
\fancyfoot[C]{\oldstylenums{\thepage}}

\smallskip

Let us consider the Lobatto quadrature \cite{gautschi1994} with two internal control points.
It requires three coefficients $\alpha$, $\beta$ and $\xi$ such that
\moneq
     \int_0^1 f(t) \dd t \approx \alpha (f(0)+f(1))+\beta (f(\xi)+f(1-\xi)).
\label{eq:lobattoQuad4points}
\monend 
If $\xi=\frac12$, Simpson's quadrature is found \cite{dubois2023_gsi,dubois2023_fvca,rojas2024_axioms}.
Let us then suppose that\br
$0<\xi<\frac{1}{2}$ so that two internal control points are satisfied. When $f(t)=1$ or $f(t)=t$,
the quadrature (\ref{eq:lobattoQuad4points}) leads to the restriction that $\alpha+\beta=\frac12$.
When $f(t)=t^2$, we have $\alpha+\beta(1-2\xi+2\xi^2)=\frac{1}{3}$, so $\beta \, \xi(1-\xi)=\frac{1}{12}$.
When $f(t)=t^4$, we have $\beta \, \xi (1-\xi)(\xi^2 - \xi +2) = \frac{3}{10}$, and by recalling the value of
$\beta \, \xi(1-\xi)$, the internal control point can be solved using $\xi(1-\xi) =\frac{1}{5}$.
The parameters $\alpha$ and $\beta$ are obtained by applying the value of $\xi$ on $\beta \,\xi(1-\xi)=\frac{1}{12}$:
$ \, \xi = \frac{1}{2}-\frac{\sqrt{5}}{10} $,
$ \, \alpha = \frac{1}{12} $, $ \, \beta = \frac{5}{12}$.
Therefore, the proposed Lobatto quadrature is
\moneq
\label{eq:lobattoQuad}
        \int_0^1 f(t) \dd t \approx \frac{1}{12} (f(0)+f(1))+\frac{5}{12}(f(\xi) + f(1-\xi)), 
        \qquad 
        \xi = \frac{1}{2}-\frac{\sqrt{5}}{10}.
\monend 
Lobatto's quadrature (\ref{eq:lobattoQuad}) is accurate for polynomials up to degree $2 \, n-3$,
where $ \, n \, $ is the number of integration points. In the present case $ \, n=4$, so the proposed Lobatto
quadrature is exact for the integration of polynomials up to degree five This can be easily verified
using the obtained values of $ \, \alpha$, $ \, \beta \, $ and $ \, \xi$.

\newpage 
\bigskip \bigskip    \noindent {\bf \large    3) \quad Cubic Interpolation} 

\smallskip 

Lagrange's $P_3$ polynomials \cite{gautschi1994} are selected for the finite-elements-based internal
interpolation of a time interval $[0,h]$, where $h$ is the step size.
For $ \, 0 \leq  \theta \leq 1$, the following four basis functions
\moneq
\left\lbrace
\begin{array}{rcl} 
\varphi_0(\theta) &=& 5  (\theta-\xi) \, (\theta-(1-\xi)) \, (1-\theta) \\
\varphi_\xi(\theta)   &=& -5 \, \sqrt{5} \, \, \theta \, (1-\theta) \, (\theta-(1-\xi))  \\
\varphi_{1-\xi}(\theta) &=& 5 \sqrt{5} \, \theta \, (1-\theta) \, (\theta-\xi) \\
\varphi_1(\theta)   &=& 5 \, \theta \, (\theta - \xi) \, (\theta-(1-\xi))
\end{array}\right.    
\label{eq:basisFunctions}
\monend 
evaluate to $1$ when the control point is the same as the index and $0$ at the other integration points.

\smallskip \noindent 
Taking $ \, t=h\theta$, where $\, h \, $ is the step size, one can build a function $t\longmapsto q(t)$ on the interval~$ \, [0,h] \, $
using the above $ \, P_3 \, $ finite elements:
\moneqstar
q(t) = q_\ell \, \varphi_0(\theta) + q_\xi \, \varphi_\xi(\theta)+ q_{1-\xi} \, \varphi_{1-\xi}(\theta) + q_r \, \varphi_1(\theta).
\monendstar 
This is a vectorial and third-order function on $[t,t+h]$. Note that $ \, q(0)=q_\ell$,
$ \, q(\frac h2-\frac{\sqrt{5}h}{10})=q_\xi$, $ \, q(\frac h2+\frac{\sqrt{5}h}{10})=q_{1-\xi}$, and $ \, q(h)=q_r$
so the basis functions (\ref{eq:basisFunctions}) are well adapted to the chosen quadrature.
The time derivative is a second-order polynomial concerning $ \, \theta$:
\moneqstar
{{\dd q}\over{\dd t}}  =\frac{1}{h} \, \left[ q_\ell  \, \varphi_0^\prime(\theta)
  + q_\xi  \, \varphi_\xi^\prime(\theta) + q_{1-\xi} \,   \varphi_{1-\xi}^\prime(\theta) + q_r  \, \varphi_r^\prime(\theta)\right].
\monendstar 
%

\bigskip \bigskip    \noindent {\bf \large    4) \quad The Harmonic Oscillator Case}


\smallskip 

Let us consider a dynamical system described by a state $q(t)$ which is a single real variable
for $ \, 0\leq  t \leq  T$. The continuous action $S_c$ introduces a Lagrangian quantity $L$ and is defined by
\moneq
	S_c = \int_0^T  L \Big( {{\dd q}\over{\dd t}} ,\, q \Big) \, \dd t.
\label{eq:continuousAction}
\monend 
For the harmonic oscillator, the Lagrangian is
$ \, L \big( {{\dd q}\over{\dd t}} ,\, q \big)=\frac m2 \, \big(  {{\dd q}\over{\dd t}} \big)^2 - V(q)$,
where the terms on the right are the kinetic and potential energies, in said order.
The action (\ref{eq:continuousAction}) is discretized by splitting the time interval $ \, [0,T] \, $ into $ \, N \, $
elements of equal length $ \, h=\frac TN$. The approximation $q_j$ of $q(t_j)$ will be given at each discrete time instance $t_j=j\,h$.

\bigskip  \monitem {\bf Discrete Lagrangian}

\smallskip 

Using the Lobatto quadrature (\ref{eq:lobattoQuad}), the kinetic energy is viewed as an integral function:  
$   \,     K_d = \frac{1}{h} \int_0^h \frac{1}{2} \, m \big(   {{\dd q}\over{\dd t}} \big)^2 \, \dd t =
\frac{1}{2} m \int_0^1 \big(\frac{1}{h}  {{\dd q}\over{\dd\theta}} \big)^2 \, \dd \theta $, 
by changing the integration variable. A symmetric quadratic form emerges as
\moneq 
K_d = \frac{m}{2h^2} \Big( \alpha \, 
\Big[ \Big({{\dd q}\over{\dd \theta}}(0) \Big)^2 + \Big({{\dd q}\over{\dd \theta}}(1) \Big)^2 \Big] 
+ \beta \, \Big[ \Big({{\dd q}\over{\dd \theta}}(\xi) \Big)^2  +   \Big({{\dd q}\over{\dd \theta}}(1-\xi) \Big)^2  \Big]  \Big) .
\label{eq:discreteKE}
\monend
Then $ \, K_d \equiv  \frac{m}{2h^2} {q_\theta}\tran K q_\theta $,
where $ \, q_\theta = \left( q_0, \,q_\xi,\, q_{1-\xi},\, q_1 \right)^{\rm t} $ and
\moneq 
K = \left(  \begin{array} {cccc}  
        \frac{13}{3} & -\frac{5\sqrt{5}}{4}-\frac{25}{12} & \frac{5\sqrt{5}}{4}-\frac{25}{12} & -\frac{1}{6} \\
        -\frac{5\sqrt{5}}{4}-\frac{25}{12} & \frac{25}{3} & -\frac{25}{6} & \frac{5\sqrt{5}}{4}-\frac{25}{12} \\
        \frac{5\sqrt{5}}{4}-\frac{25}{12} & -\frac{25}{6} & \frac{25}{3} & -\frac{5\sqrt{5}}{4}-\frac{25}{12} \\
        -\frac{1}{6} & \frac{5\sqrt{5}}{4}-\frac{25}{12} & -\frac{5\sqrt{5}}{4}-\frac{25}{12} & \frac{13}{3}
\end{array} \right)  \,.
\label{eq:rigidityKE}
\monend 

\noindent 
When the potential energy function $q \longmapsto V(q)$ is a polynomial of degree
$ \, \leq  5$, the quadrature (\ref{eq:lobattoQuad}) integrates it accurately and
\moneq 
        U \equiv \frac{1}{h} \int_0^h V(q(t)) \, \dd t = \alpha [V(q_0)+V(q_1)]+\beta[V(q_\xi)+V(q_{1-\xi})].
\label{eq:discretePE}
\monend 
In this example, the potential energy function is $V(q) = \frac{1}{2} m \omega^2 q^2$. The discrete Lagrangian
$ \, L_d \equiv h(K_d - U)$, can be expressed as a function of the degrees of freedom within the interval
$ \, [0,h]$: $ \, q_\ell = q(0)$, $ \, q_\xi=q(\xi h)$, $ \, q_{1-\xi} = q((1-\xi)h)$, and
$ \, q_r=q(h)$. The Lagrangian   $ \, L_d \equiv  L_d\left(q_\ell,q_\xi,q_{1-\xi},q_r \right) \, $ is given by
\moneq
L_d = \left\{  \begin{array} {l}
\frac{m}{12 \, h} \Big[ 26 \, \left({q_\ell}^2+{q_r}^2\right) - 2 \, q_\ell q_r + 50\left({q_\xi}^2-q_\xi q_{1-\xi}
  + {q_{1-\xi}}^2 \right) \\ 
\qquad   -25 \, (q_\ell+ q_r)(q_\xi + q_{1-\xi}) - 15\sqrt{5} (q_\ell- q_r)(q_\xi - q_{1-\xi}) \Big] \\ 
\qquad \qquad -\frac{m\,h\omega^2}{24} \, \Big[ {q_\ell}^2 + {q_r}^2 + 5\left( {q_\xi}^2 + {q_{1-\xi}}^2 \right) \Big].
\end{array} \right.  
\label{eq:discreteLagrangian}
\monend
The internal degrees of freedom $q_\xi$ {and} $q_{1-\xi}$ are eliminated by solving the system formed by
the discrete Euler-Lagrange equations inside the interval, that is
\moneqstar 
    \frac{\partial L_d}{\partial q_\xi}=0 ,\qquad \frac{\partial L_d}{\partial q_{1-\xi}}=0.
\monendstar
The explicit expressions of $q_\xi$ and $q_{1-\xi}$ are  
\moneq
\left\{  \begin{array} {rcl}
  q_\xi & = & \frac{1}{\delta} \left[ -5\left( h^2 \omega ^2 - 30\right)(q_r+q_\ell)
    + 3 \sqrt{5} \left(h^2 \omega ^2 -10\right)(q_r-q_\ell) \right] 	\\  [0.5ex]
  q_{1-\xi}  & = &
\frac{1}{\delta} \left[ -5 \left(h^2 \omega ^2 -30\right)(q_r+q_\ell)
  - 3 \sqrt{5} \left(h^2 \omega ^2 -10\right)(q_r-q_\ell) \right] \,, 
\end{array} \right.  
\label{eq:qxiq1mxi}
\monend
with a determinant
\moneq
	\delta = (h^2\omega^2-30) \, (h^2\omega^2-10) > 0 \quad {\rm{for}} \quad 0<h\omega < \sqrt{10} \,,
	\label{eq:stabilityCondition}
\monend 
{indicating} conditional stability. By using the expressions of equations (\ref{eq:qxiq1mxi}) within the discrete
Lagrangian (\ref{eq:discreteLagrangian}), a reduced Lagrangian $ \, L_r(q_\ell,q_r) \, $ is revealed:
\moneqstar
L_r =  \left\{  \begin{array} {l}
\frac{1}{6 \, h\,\delta}
\Big[ \frac m4 \, \big(-h^6 \, \omega^6 + 92 \, h^4\omega^4 - 1680 \, h^2 \, \omega^2 + 3600 \big) \,
  \big( {q_\ell}^2  + {q_r}^2 \big)  \\
\qquad \qquad - m \, \big( \, h^4 \, \omega^4 + 60\, h^2\, \omega^2 + 1800 \big) \, q_\ell \, q_r \Big]. 
\end{array} \right.
\monendstar

%

\bigskip  \monitem {\bf Euler-Lagrange Equations}

\smallskip 

A discrete form $S_d$ of the continuous action $S_c$ (\ref{eq:continuousAction}) involves the reduced Lagrangian:
\moneqstar 
S_d = \sum_j L_r(q_j,q_{j+1})  \,.
\monendstar 
It is stationary with respect to the internal degree of freedom $ \, q_j \, $
when $  \, {{\partial S_d}\over{\partial q_j}} =0$. This leads to the discrete Euler-Lagrange equations 
\moneq 
{{\partial}\over{\partial q_r}} {L_r(q_{j-1},q_j)} + {{\partial}\over{\partial q_\ell}} {L_r(q_{j},q_{j+1})} = 0 \,,
\label{eq:centerEL}
\monend 
which become
\moneq \left\{  \begin{array} {l}
\frac{1}{h^2} \, \left( q_{j-1} - 2q_j + q_{j+1}\right) 
	+ \frac{\omega^2}{30}\left( q_{j-1} + 28q_j + q_{j+1} \right)   \\ [.5ex] 
\qquad  + \frac{\omega^4 \, h^2}{1800} \, \left( q_{j-1} - 92q_j + q_{j+1} \right) + \frac{\omega^6 \, h^4}{1800} q_j = 0.
\end{array} \right.
\label{eq:discreteCenterEL}
\monend 
The truncation order is obtained by replacing the discrete variables $ \, q_{j+1}$, $ \,, q_j$ and $ \, q_{j-1}$
by the solution of the differential equation at $ \, t_{j}+h$, $ \, t_j$ and $ \, t_j-h$. Then
\moneqstar 
\left\{ \begin{array} {rcl}
  q_{j+1} & = & q_j + h \, {{\dd q}\over{\dd t}} + \frac{h^2}{2} \, {{\dd^2 q}\over{\dd t^2}}
  + \frac{h^3}{6} \, {{\dd^3 q}\over{\dd t^3}} + \frac{h^4}{24} \, {{\dd^4 q}\over{\dd t^4}}
  + \frac{h^5}{120} \, {{\dd^5 q}\over{\dd t^5}}  + \frac{h^6}{720}  \, {{\dd^6 q}\over{\dd t^6}} + {\rm O} (h^7)  \\ [0.5ex]
  q_{j-1} & = & q_j - h \, {{\dd q}\over{\dd t}} + \frac{h^2}{2} \, {{\dd^2 q}\over{\dd t^2}}
  - \frac{h^3}{6} \, {{\dd^3 q}\over{\dd t^3}} + \frac{h^4}{24} \, {{\dd^4 q}\over{\dd t^4}}
  - \frac{h^5}{120} \, {{\dd^5 q}\over{\dd t^5}}  + \frac{h^6}{720}  \, {{\dd^6 q}\over{\dd t^6}} + {\rm O} (h^7) \, . 
\end{array} \right.
\monendstar 
Under these conditions, the left-hand side of equation (\ref{eq:discreteCenterEL}) does not vanish but defines
the truncation error $\mathcal{T}_h(q_j)$. 
The scheme (\ref{eq:discreteCenterEL}) is sixth-order accurate on truncation error:
$ \, {\mathcal{T}}_h(q_j) = - \frac{1}{21600} \, \omega^8 \, h^6 \, q_j + {\rm O} (h^8) $.
%

\smallskip \noindent 
Regarding numerical stability, a condition is obtained from solving a characteristic polynomial of equation (\ref{eq:discreteCenterEL}),
\moneqstar 
\frac{1}{h^2} \, \left( 1 - 2r + r^2\right) + \frac{\omega^2}{30} \, \left( 1 + 28r + r^2 \right) +
\frac{\omega^4 h^2}{1800} \, \left( 1 - 92r + r^2 \right) + \frac{\omega^6 h^4}{1800} r = 0 \,.
\monendstar 
The corresponding discriminant is
\moneqstar 
\Delta = \frac{\omega^2}{h^2} \, \left(h^2 \omega ^2-10\right) \, \left(h^2 \omega^2-30\right) 						
\, \left(h^2 \omega ^2-60\right) \, \left(h^4 \omega ^4-84 h^2 \omega^2+720\right) \, .
\monendstar 
The polynomial $ \, \left(h^4 \omega ^4-84 h^2 \omega^2+720\right) \, $ has two real roots in $h^2\omega^2$. Therefore,
\moneq
\Delta < 0 \quad {\rm when} \quad 0 < h \, \omega < \sqrt{42-6\sqrt{29}} \, . 
\label{eq:stabilityCondition2}
\monend
This restriction supersedes the previous stability condition (\ref{eq:stabilityCondition}).
Under this new restriction, two complex conjugate roots of unit module are obtained, guaranteeing
numerical stability. Let us remark that
$ \, \sqrt{42-6\sqrt{29}} \approx \pi$. As such, $ \, h<\frac T2 \, $ when $ \, \omega = \frac{2\pi}{T}$.
The method remains stable by using a little more than two points per period of oscillation in this case.

\bigskip  \monitem {\bf  Symplectic Structure}

\smallskip

The generalized momentum is defined on the right  by
\moneq 
\label{eq:rightMomentum}
	p_r = \diffp{L_r}{q_r}\left( q_\ell, q_r\right).
\monend
In the case of the harmonic oscillator,
\moneqstar
p_r = m \frac{q_r-q_\ell}{h} + \frac{m}{6 \, \delta} \,
\Big[ -300 \, h \, \omega^2 \, \left(q_\ell+2\,q_r\right) + 5 \, h^3 \, \omega^4 \left(q_\ell+8\, q_r\right)
  -\frac12 \, h^5 \, \omega^6 q_r \Big].
\monendstar 
By noticing that $p_j = \diffp*{L_r\left(q_{j-1},q_j\right)}{q_r}$, equation (\ref{eq:centerEL}) gives 
$ \, p_j = -\diffp*{L_r\left(q_j,q_{j+1}\right)}{q_\ell} $, 
and $ \, p_{j+1} \, $ is calculated according to (\ref{eq:rightMomentum}).
A discrete system involving the momentum and the state is obtained:
\moneq
\label{eq:discreteSystem}
\displaystyle \left\{  \begin{array} {rcl}
\displaystyle  \frac{p_{j+1}-p_j}{h} &=& \displaystyle 
-m \, \omega ^2 \, \frac{\left(h^2 \, \omega^2 -60 \right)}{12 \left(h^2 \, \omega^2-10\right)} \, \left(q_j+q_{j+1}\right), \\ [.6ex]
 \displaystyle  \frac{q_{j+1}-q_j}{h} &=& \displaystyle
\frac{24 \, (30 - h^2\, \omega^2)}{\left(h^4 \, \omega ^4-84 h^2 \, \omega ^2+720\right)}  \Big(\frac{p_j+p_{j+1}}{2\,m} \Big).
\end{array} \right.
\monend
The system (\ref{eq:discreteSystem})  gives the recurrence iteration scheme 
\begin{equation}
 \left( \begin{array} {c} p_{j+1}  \\ q_{j+1} \end{array} \right) 
= \mathbf{\Phi} \, 
 \left( \begin{array} {c} p_j \\ q_{j} \end{array} \right) \,,\quad 
\mathbf{\Phi} = \frac{1}{\widetilde{\delta}}  \, \left( \begin{array} {cc} a & b \\ c & a  \end{array} \right) \,, 
\label{eq:symplecticSystem}
\end{equation}
with
\moneq 	\left\{ \begin{array}{rcl}
\widetilde{\delta} & = & 1+ \frac{h^2\omega^2}{30} + \frac{h^4\omega^4}{1800} \\ [0.5ex]
a & = & 1 -\frac{7}{15}h^2\omega^2 +\frac{23}{900}h^4\omega^4 -\frac{1}{3600}h^6\omega^6 \\  [0.5ex]
b & = & \frac{m\, h\omega^2}{43200} \left(h^2\omega^2 -60\right)\left(h^4\omega^4 -84h^2\omega^2 +720\right) \\  [0.5ex]
c & = & \frac hm \left( 1-\frac2{15}h^2\omega^2 +\frac1{300}h^4\omega^4 \right).
\end{array} \right.
\label{eq:symplecticCoeffs}
\end{equation}
Let us remark that $\det \mathbf{\Phi} = 1$, so the discrete flow (\ref{eq:symplecticSystem}) is symplectic
(see \textit{e.\,g.} \cite{sanzSerna1992}).

\noindent 
The harmonic oscillator preserves the energy quantity $ \, H(p,q)= \frac{1}{2m}p^2 + \frac{m\omega^2}{2}q^2 $.
According to \cite{dubois2023_gsi}, and using (\ref{eq:symplecticCoeffs}) the Lobatto scheme preserves the discrete energy 
\moneq  
	H_d(p,q) = \frac{c}{2\,\widetilde{\delta}} \, p^2 - \frac{b}{2 \, \widetilde{\delta}} \, q^2.
\label{eq:discreteEnergy}
\monend

\smallskip \noindent 
Under the stability condition (\ref{eq:stabilityCondition2}), we have the inequality $ \, -b>0$.

\bigskip \bigskip    \noindent {\bf \large    5) \quad The Nonlinear Pendulum Case}

\vskip -0.2 cm\smallskip

Let us now consider a nonlinear pendulum system described by a the state $q(t)$
and a constant mass $m$. The continuous action (\ref{eq:continuousAction}) and Lagrangian
coincide with that of the harmonic oscillator of the previous section. However,
the potential $V(q)=m\omega^2(1-\cos q)$ in this case is nonlinear. The second order
dynamics are given by $\diff[2] qt + \omega^2 \sin q = 0$. An analytical solution of
this problem involves special elliptic functions and is established in, \emph{e.\,g.}, \cite{dubois2023_fvca}.

\bigskip  \monitem {\bf Discrete Lagrangian and Euler-Lagrange Equations}

\smallskip

As before, the discrete Lagrangian
$ \, L_d \equiv h(K_d - U)$ using equations (\ref{eq:discreteKE}), (\ref{eq:rigidityKE}) and (\ref{eq:discretePE}); can be expressed as a function of the degrees of freedom within the interval
$ \, [0,h]$: $ \, q_\ell = q(0)$, $ \, q_\xi=q(\xi h)$, $ \, q_{1-\xi} = q((1-\xi)h)$, and
$ \, q_r=q(h)$. For this example, $ \, L_d \equiv  L_d\left(q_\ell,q_\xi,q_{1-\xi},q_r \right) \, $ is given by
\moneq
L_d = \left\{  \begin{array} {l}
\frac{m}{12 \, h} \Big[ 26 \, \left({q_\ell}^2+{q_r}^2\right) - 2 \, q_\ell q_r + 50\left({q_\xi}^2-q_\xi q_{1-\xi}
  + {q_{1-\xi}}^2 \right) \\ 
\qquad   -25 \, (q_\ell+ q_r)(q_\xi + q_{1-\xi}) - 15\sqrt{5} (q_\ell- q_r)(q_\xi - q_{1-\xi}) \Big] \\ 
\qquad \qquad -\frac{h}{12} \, \Big[ V_\ell + 5\left( V_\xi + V_{1-\xi} \right) + V_r  \Big].
\end{array} \right. ,  
\label{eq:discreteLagrangianNL}
\monend
where $V_\ell\equiv V(q_\ell)$, $V_\xi\equiv V(q_\xi)$, $V_{1-\xi}\equiv V(q_{1-\xi})$, and $V_r\equiv V(q_r)$. Unlike the linear, harmonic oscillator case, internal variables $\{q_{1-\xi},q_\xi\}$ cannot be eliminated.

Euler-Lagrange equations, coming from the stationary action $\delta S_d = 0$, are first established for arbitrary variations $\delta q_\xi$ and $\delta q_{1-\xi}$
in the interval $[t_j,t_{j+1}]$: 
\moneq
	{{\partial L_d}\over{\partial q_\xi}}=0; \qquad {{\partial L_d}\over{\partial q_{1-\xi}}}=0.
	\label{eq:ELinternalNL}
\monend
Using the expression (\ref{eq:discreteLagrangianNL}) of the discrete Lagrangian on equation (\ref{eq:ELinternalNL}),
the sums

\noindent $2 \, {{\partial L_d}\over{\partial q_\xi}} + {{\partial L_d}\over{\partial q_{1-\xi}}} = 0 \, $ and
$ \,{{\partial L_d}\over{\partial q_\xi}} +  2 \, {{\partial L_d}\over{\partial q_{1-\xi}}} = 0 $, respectively give
\moneq
	\left\{
	\begin{array}{rcl}
		q_\xi - \frac{h^2}{30m}\left( V_{1-\xi}' + 2 V_\xi'\right) 
		& = 
		& \frac{5+\sqrt{5}}{10}q_\ell + \frac{5-\sqrt{5}}{10}q_r,
		\\
		q_{1-\xi} - \frac{h^2}{30m}\left( 2V_{1-\xi}' +  V_\xi'\right) 
		& = 
		& \frac{5-\sqrt{5}}{10}q_\ell + \frac{5+\sqrt{5}}{10}q_r,
	\end{array}
	\right.
	\label{eq:internesNL}
\monend
where we adopt the notation $V_x'\equiv \diffp{V}{q}[q_x]$. The expressions above implicitly define $\{q_\xi,q_{1-\xi}\}$ within the interval as a function of the border values $\{q_\ell,q_r\}$.

On the borders, the Euler-Lagrange equations are
\moneq
	\diffp{L_d}{q_r}\!\left(q_{j-1},q_{j-1+\xi},q_{j-\xi},q_j \right)
	 + \diffp{L_d}{q_\ell}\!\left(q_j,q_{j+\xi},q_{j+1-\xi},q_{j+1} \right)
	 = 0.
	\label{eq:ELequationsNL}
\monend
By injecting the expressions of $\{q_{1-\xi},q_\xi\}$ from (\ref{eq:internesNL}), the Euler-Lagrange equations for the borders read
\moneq
	\begin{array}{l}
	\frac{1}{h^2}\left( q_{j-1}-2q_j+q_{j+1} \right)
	+ \frac{1}{24m}\Big[
		\left( 5-\sqrt{5} \right) V_{j-1+\xi}'
		+ \left( 5+\sqrt{5} \right) V_{j-\xi}'
		+ 4 V_{j}'
		\\
		\qquad \qquad \qquad \qquad \qquad \qquad \quad \
		+ \left( 5+\sqrt{5} \right) V_{j+\xi}'
		+ \left( 5-\sqrt{5} \right) V_{j+1-\xi}' \Big] = 0.
	\end{array}
\monend
Momenta are given on the right by $p_r = \diffp{L_d}{q_r}\!\left(q_\ell,q_\xi,q_{1-\xi},q_r \right)$.
As such,
\moneqstar 
p_{j+1} = \diffp{L_d}{q_r}\!\left(q_j,q_{j+\xi},q_{j+1-\xi},q_{j+1} \right) \, . 
\monendstar 
Then, from the Euler-Lagrange equations (\ref{eq:ELequationsNL}), one can identify
\moneqstar 
p_j = -\diffp{L_d}{q_\ell}\!\left(q_j,q_{j+\xi},q_{j+1-\xi},q_{j+1} \right) \, .
\monendstar 
From these relationships, the discrete Hamiltonian dynamics can be deduced:
\moneq
	\arraycolsep=0.75pt
	\begin{array}{rcl}
		p_{j+1}-p_j 
		+ \frac{h}{12}\! \left[ 
		V_j' + 5\left( V_{j+\xi}' + V_{j+1-\xi}'\right) + V_{j+1}' 
		\right]
		& =
		& 0,
		\\
		q_{j+1}-q_j 
		- \frac{h^2}{24m}\! \left[ 
		V_{j+1}' + \sqrt{5}\!\left( V_{j+1-\xi}' - V_{j+\xi}'\right) - V_{j}'
		\right]
		-\frac{h}{2m}\!\left(p_{j+1} + p_j \right)
		& =
		& 0.
	\end{array}
	\label{eq:discreteDynamics}
\monend
The system (\ref{eq:internesNL}), (\ref{eq:discreteDynamics}) is written under the form $F_L\!\left(q_{j+\xi},q_{j+1-\xi},p_{j+1},q_{j+1}\right)=0$. We have
\moneqstar
	\arraycolsep=1ex
	\dd F_L\!\left(q_{\xi},q_{1-\xi},p,q\right) = 
	\left(
	\begin{array}{cccc}
		1-\frac{h^2}{15m}V_\xi''
		& -\frac{h^2}{30m}V_{1-\xi}''
		& 0
		& \frac{-5+\sqrt{5}}{10}
		\\
		-\frac{h^2}{30m}V_{1-\xi}''
		& 1-\frac{h^2}{15m}V_{1-\xi}''
		& 0
		& -\frac{5+\sqrt{5}}{10}
		\\
		\frac{5h}{12}V_\xi''
		& \frac{5h}{12}V_{1-\xi}''
		& 1 
		& \frac{h}{12}V''
		\\
		\frac{\sqrt{5}h^2}{24m}V_\xi''
		& -\frac{\sqrt{5}h^2}{24m}V_{1-\xi}''
		& -\frac{h}{2m}
		& 1-\frac{h^2}{24m}V''
	\end{array}
	\right)
\monendstar
where we adopt the notation $V_x''\equiv \diffp[2]Vq[q_x]$ and $V''\equiv V''\!(q)$.
Let us remark that

\noindent $\det \, (\!\dd F_L) = 1+\frac{h^2}{60m}\left( V_\xi'' + V_{1-\xi}'' \right) + \frac{h^4}{1800m}V_\xi''V_{1-\xi}'' \, $
tends to $1$ as $h$ diminishes. After calculating the inverse of the above Jacobian matrix,
the four-equation nonlinear system
\moneqstar 
F_L\!\left(q_{j+\xi},q_{j+1-\xi},p_{j+1},q_{j+1}\right)=0
\monendstar
is solved with Newton's algorithm.
We have observed machine precision convergence at the fifth iteration with the proposed scheme.

\bigskip  \monitem {\bf Symplectic Structure}

\smallskip
From equations (\ref{eq:internesNL}), we have
\moneqstar
	\left\{
	\begin{array}{rcl}
		\delta q_{j+\xi}
		& = 
		& \frac{1}{\Delta}\Big[ 
			30m^2\left[ 
				(5+\sqrt{5})\delta q_j 
				- (-5+\sqrt{5})\delta q_{j+1}
			\right] 
		\\
		&
		&	\qquad - mh^2 V_{j+1-\xi}'' \left[ 
				(5+3\sqrt{5})\delta q_j 
				+ (5-3\sqrt{5})\delta q_{j+1}
			\right]
		  \Big],
		\\
		\delta q_{j+1-\xi}
		& = 
		& \frac{1}{\Delta}\Big[ 
			-30m^2\left[ 
				(-5+\sqrt{5})\delta q_j 
				- (5+\sqrt{5})\delta q_{j+1}
			\right] 
		\\
		&
		&	\qquad - mh^2 V_{j+\xi}'' \left[ 
				(5-3\sqrt{5})\delta q_j 
				+ (5+3\sqrt{5})\delta q_{j+1}
			\right]
		  \Big],
	\end{array}
	\right.
\monendstar
where $\Delta = 300m^2-20mh^2 \left( V_{j+\xi}'' + V_{j+1-\xi}'' \right) + h^4V_{j+\xi}''V_{j+1-\xi}''$. Then, by differentiating equations (\ref{eq:discreteDynamics}) and using $\{\delta q_{j+\xi},\delta q_{j+1-\xi}\}$ above, we obtain a system that can be expressed as
\begin{equation}
 A \left( \begin{array} {c} \delta p  \\ \delta q \end{array} \right)_{j+1} 
= B 
 \left( \begin{array} {c} \delta p \\ \delta q \end{array} \right)_j \,
\label{eq:symplecticSystemNL}
\end{equation}
where $\det A = \det B = \frac{1}{6\Delta}\left[1800m^2 + 30mh^2 \left( V_{j+\xi}'' + V_{j+1-\xi}'' \right) + h^4V_{j+\xi}''V_{j+1-\xi}''\right]$, and the relationship
\moneqstar
	\diffp{p_{j+1}}{p_j} \diffp{q_{j+1}}{q_j} - \diffp{p_{j+1}}{q_j} \diffp{q_{j+1}}{p_j} = 1
\monendstar
is established. The Lobatto scheme (\ref{eq:discreteDynamics}) is symplectic.

\bigskip \bigskip    \noindent {\bf \large    6) \quad Numerical Experiments}

\smallskip 

The harmonic oscillator evolution was simulated using the Lobatto scheme (\ref{eq:symplecticSystem}) (\ref{eq:symplecticCoeffs}).
Results for both $ \, N = 3$ meshes are displayed on Fig. \ref{fig:simulations}, compared with the exact
solutions $ \, q(t)=\frac{\pi}{2} \, \cos (\omega t) \, $ and $ \, p(t) = -\frac{\pi}{2} \, m \, \omega \, \sin(\omega t) \, $
on a period $ \, T=1 $.
Comparable results from both the implicit midpoint 
and Simpson schemes are also provided for reference (see \textit{e.\,g.} \cite{dubois2023_gsi}).
Quantitative errors with the $ \, \ell^\infty $ norm are given in Tab. \ref{tab:errorNorms}.
An asymptotic order of convergence of 6 is estimated for the momentum, states, and energy.
Fig. \ref{fig:errorEvolution} shows the $\ell^\infty$ energy error norm evolution across $ \num{1e5} $  periods of motion.

\begin{figure}[H] 
\centering
\includegraphics[]{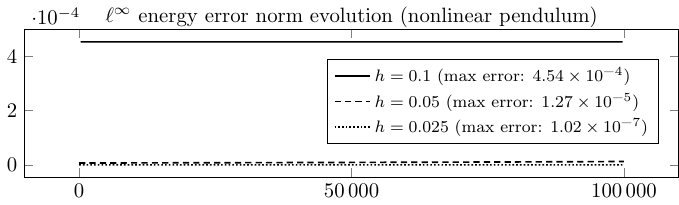}
\setlength{\abovecaptionskip}{0ex}
\caption{Over $ \num{5e3}$ periods, the $\ell^\infty$ energy error growth rate is of: 
\num{7.99e-13}   when $ \, h=0.1$,  \num{5.39e-11} when $ \, h=0.05$, and  \num{5.33e-15} when $ \, h=0.025$.}%
	\label{fig:errorEvolutionNL}%
\end{figure}

{\noindent The nonlinear pendulum evolution was simulated using the nonlinear Lobatto scheme (\ref{eq:internesNL}) (\ref{eq:discreteDynamics}).
Results for $\, N = 5$ meshes are displayed on Fig. \ref{fig:simulationsNL}, compared with the exact
solution involving special functions (see \emph{e.\,g.} \cite{dubois2023_fvca}).
Comparable results from both the implicit midpoint and Simpson schemes are also provided for reference (see \textit{e.\,g.} \cite{dubois2023_fvca}).
Quantitative errors with the $ \, \ell^\infty $ norm are given in Tab. \ref{tab:errorNormsNL}.
An asymptotic order of convergence of 6 is estimated for the momentum, states, and energy.
Fig. \ref{fig:errorEvolutionNL} shows the $\ell^\infty$ energy error norm evolution across $ \num{1e5} $  periods of motion.

\smallskip 
\begin{figure}[H] 
\centering
\includegraphics[width=\columnwidth]{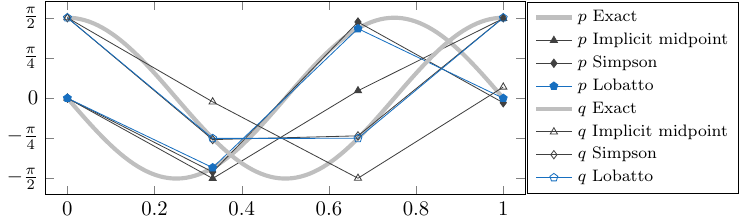}


%
\setlength{\abovecaptionskip}{0.5ex}
\caption{Harmonic oscillator evolution for the momentum $ \, p \, $ and state $ \, q $.
Comparison of the Lobatto symplectic scheme against the exact solution for $N = 3$ meshes.
Lobatto's solutions are very close to the exact ones. 
Momentum data have been rescaled.}%
\label{fig:simulations}%
\end{figure}

\begin{table}[H]
	\setlength{\belowcaptionskip}{1ex}
	\caption{Errors in the maximum norm. The Lobatto approximation is sixth-order accurate and preserves
          a discrete energy. Comparable values from both the implicit midpoint and Simpson schemes can be found in \cite{dubois2023_gsi}.
          The estimated convergence order is the closest integer $\alpha$ measuring the ratio of successive errors in a given line
          by a negative power of 2 of the type $2^{-\alpha}$.}
	\label{tab:errorNorms}
	\centering
		\sisetup{
					scientific-notation = true,
					round-mode = places,
					round-precision = 2,
					table-alignment-mode = marker,
					table-format =1.2e-2,
					table-number-alignment = center}
		\begin{tabular}{
			l
			S
			S
			S
			c
		}
			\toprule
			Number of meshes
			& \text{10}
			& \text{20}
			& \text{40}
			& order
			\\
			\midrule
			 Momentum $p$
			& 8.95248e-6
			& 1.39279e-7
			& 2.17034e-9
			& 6
			\\
			 State $q$
			& 7.64034e-7
			& 1.19368e-8
			& 1.87641e-10
			& 6
			\\ 
			 Energy $H(p,q)$
			& 0.0000661948
			& 1.09831e-6
			& 1.69917e-8
			& 6
			\\
			 Discrete energy $H_d(p,q)$ 
			& 2.66454e-15
			& 1.33227e-15
			& 2.22045e-15
			& exact
			\\
			\bottomrule
		\end{tabular}
	\end{table}

\begin{figure}[H] 
\centering
\includegraphics[]{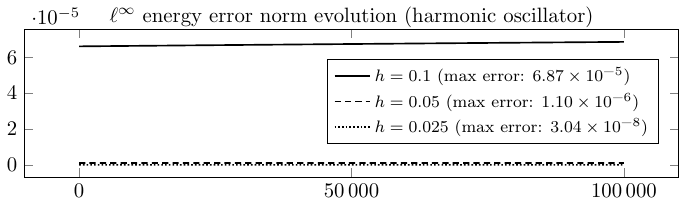}
\setlength{\abovecaptionskip}{0ex}
\caption{Over $ \num{1e5}$ periods, the $\ell^\infty$ energy error growth rate is of: 
\num{2.50e-11}   when $ \, h=0.1$,  \num{5.45e-14} when $ \, h=0.05$, and  \num{1.35e-13} when $ \, h=0.025$.}%
	\label{fig:errorEvolution}%
\end{figure}

\begin{figure}[H] 
\centering
\includegraphics[width=\columnwidth]{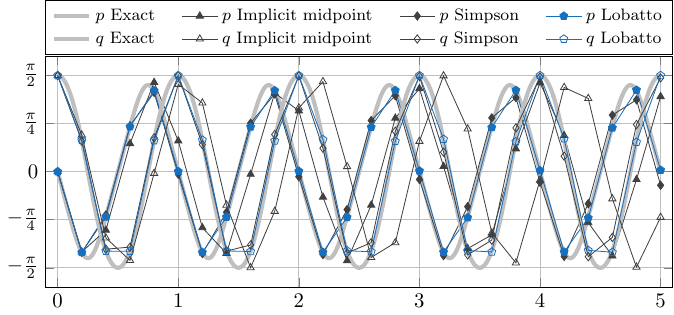}

\setlength{\abovecaptionskip}{0.5ex}
\caption{Nonlinear pendulum evolution for the momentum $ \, p \, $ and state $ \, q $.
Comparison of the Lobatto symplectic scheme against the exact solution for $N = 5$ meshes over 5 periods.
Lobatto's solutions remain very close to the exact ones.
Note that the momentum data have been rescaled.}%
\label{fig:simulationsNL}%
\end{figure}

\begin{table}[H]
	\setlength{\belowcaptionskip}{1ex}
	\caption{Errors in the maximum norm for the nonlinear pendulum simulation. The Lobatto integrator remains a sixth-order method in the nonlinear case. Comparable values for both the implicit midpoint and Simpson schemes can be found in \cite{dubois2023_fvca}.
	}
	\label{tab:errorNormsNL}
	\centering
		\sisetup{
					scientific-notation = true,
					round-mode = places,
					round-precision = 2,
					table-alignment-mode = marker,
					table-format =1.2e-2,
					table-number-alignment = center}
		\begin{tabular}{
			l
			S
			S
			S
			c
		}
			\toprule
			Number of meshes
			& \text{50}
			& \text{100}
			& \text{200}
			& Order
			\\
			\midrule
			 Momentum $p$
			& 2.83174e-9
			& 4.56684e-11
			& 7.0699e-13
			& 6
			\\
			 State $q$
			& 4.21838e-10
			& 6.69165e-12
			& 1.05693e-13
			& 6
			\\ 
			 Energy $H(p,q)$
			& 6.23384e-10
			& 1.02788e-11
			& 1.58925e-13
			& 6
			\\
			\bottomrule
		\end{tabular}
	\end{table}
}

\bigskip \bigskip    \noindent {\bf \large    5) \quad Conclusions}

  We have constructed a variational integrator based on Lobatto's quadrature and applied it to the
  least action principle. Using this method, the discrete formulation of two fundamental problems is given:
  the harmonic oscillator and the nonlinear pendulum for which the analytical solution\textcolor{magenta}{s are} known.
After completing this work, the authors noticed that {the linear harmonic oscillator} scheme (\ref{eq:symplecticSystem}) (\ref{eq:symplecticCoeffs})
was published in \cite{oberBlobaum2015}. Our analysis coincides with that of \cite{oberBlobaum2015}
{but presents different simulation results and explicits a preserved quantity. Additionally,
  the nonlinear integrator is provided by the scheme (\ref{eq:internesNL}) (\ref{eq:discreteDynamics}),
  and tested on the nonlinear pendulum example. The method is symplectic and conditionally stable.
  It uses cubic Lagrange polynomials for its internal interpolation and results to be sixth-order
  accurate for the state, the momentum, and the system energy (twice the order of the selected polynomials).
  This superconvergence is a characteristic of Galerkin variational integrators as remarked
  in \cite{oberBlobaum2015,oberBlobaum2021}. The formulation of the integrator on multi-degrees
  of freedom nonlinear systems is left for future work.}

\bigskip \bigskip      \noindent {\bf  \large  References }



\end{document}